\documentclass[12pt]{article}

\usepackage{amssymb}
\usepackage{latexsym}
\usepackage{amsthm}
\usepackage{amscd}
\usepackage{amsmath}
\usepackage{diagrams}

\newtheorem{thm}{Theorem}[section]

\newtheorem{defn-lem}[thm]{Definition-Lemma}
\newtheorem{conj}[thm]{Conjecture}

\theoremstyle{remark}
\newtheorem{rem}[thm]{Remark}

\theoremstyle{definition}
\newtheorem{defn}[thm]{Definition}
\numberwithin{equation}{section}

\allowdisplaybreaks[4]

\def \Q{{\mathbb Q}}
\def \N{{\mathbb N}}
\def \C{{\mathbb C}}

\def \Z{{\mathbb Z}}
\def \R{{\mathbb R}}
\def \P{\mathbb P}

\def\map#1.#2.{#1 \longrightarrow #2}
\def\rmap#1.#2.{#1 \dasharrow #2}

\DeclareMathOperator{\Hom}{Hom}
\DeclareMathOperator{\Hol}{Hol}

\def\fb#1.{\underset #1 \to \times}
\def\pr#1.{\Bbb P^{#1}}
\def\ring#1.{\mathcal O_{#1}}
\def\mlist#1.#2.{{#1}_1,{#1}_2,\dots,{#1}_{#2}}

\def\Hom{\operatorname{Hom}}

\def\uloopr#1{\ar@'{@+{[0,0]+(-4,5)} @+{[0,0]+(0,10)}
@+{[0,0]+(4,5)}}
  ^{#1}}
\def\dloopr#1{\ar@'{@+{[0,0]+(-4,-5)} @+{[0,0]+(0,-10)}
@+{[0,0]+(4,-5)}}
  _{#1}}

\def\rloopd#1{\ar@'{@+{[0,0]+(5,4)} @+{[0,0]+(10,0)}
@+{[0,0]+(5,-4)}}
  ^{#1}}
\def\lloopd#1{\ar@'{@+{[0,0]+(-5,4)} @+{[0,0]+(-10,0)}
@+{[0,0]+(-5,-4)}}
  _{#1}}

\long\def\ignore#1{}
\long\def\ignore#1{#1}

\begin{document}

\begin{center}

{\bf Rational homotopy stability for the spaces of rational maps}

\end{center}
\begin{center}
{Jiayuan Lin}
\end{center}
{\bf {\small abstract}} {\small Let $\Hol_{x_0}^{{\bf n}} (\C\P^1, X)$ be the space of based holomorphic maps of degree ${\bf n}$ from $\C\P^1$ into a simply connected algebraic variety $X$.  Under some condition we prove that the map $\map \Hol_{x_0}^{{\bf n}} (\C\P^1, X). \Hol_{x_0}^{d{\bf n}} (\C\P^1, X).$ obtained by compositing $f \in \Hol_{x_0}^{{\bf n}} (\C\P^1, X)$ with $g(z)=z^d, z \in \C\P^1$ induces rational homotopy equivalence up to some dimension, which tends to infinity as the degree grows.} 
\section{Introduction}
\medskip

    Rational curves play a more and more essential role either in defining quantum cohomology or in studying higher dimensional algebraic varieties (see Koll\'ar $[16]$ for example). The space of rational maps from $\C\P^1$ into certain target space $X$ has attracted a lot of attention from several mathematical branches.

    Let $\Hol_{x_0}^{n} (\C\P^1, X)$ and $F_{x_0}^{n} (\C\P^1, X)$ be the space of based holomorphic maps and the space of based continuous maps of degree $n$ from $\C\P^1$ into an algebraic variety $X$ respectively.  In $1979$, Segal $[23]$ proved that the inclusion map $\Hol_{x_0}^{n} (\C\P^1, \C\P^m) \subset F_{x_0}^{n} (\C\P^1, \C\P^m)$ is a homotopy equivalence up to dimension $n(2m-1)$. His work was motivated by the observation of Atiyah that the only critical points of an
\lq\lq energy" functional on the space of continuous maps is the space of rational maps
(where the functional achieves an absolute minimum) and by the extrapolation of finite
dimensional Morse theory to the infinite dimensional case. 
Segal's theorem was then generalized to the case $X$ a complex Grassmannian by Kirwan $[15]$, and certain $SL(n,\C)$ flags by Guest $[8]$. Later Mann and Milgram $[18]$ increased the range of the isomorphisms obtained by Kirwan for Grassmannian and treated all $SL(n,\C)$ flag manifolds in $[19]$. Similar stability theorems were proved for any generalized flag manifold $G/P$ by Boyer, Hurtubise, Mann and Milgram $[2]$ and Hurtubise $[13]$, following on a stable result of Gravesen $[7]$ and for toric varieties by Guest $[9]$. We refer to $[14]$ for a survey.

More recently, Boyer, Hurtubise and Milgram $[1]$ questioned what is the most general complex target space $X$ which admits stability theorem, which, loosely speaking, say that the homotopy (homology) groups of the space of holomorphic maps is isomorphic to the homotopy (homology) groups of the entire mapping space through a range that grows with degree, as degree moves to infinity. They noticed that there is sort of link between the stability and some rationality property of the manifold. The stability theorem should hold for certain subclass of rationally connected varieties (a variety is called rationally connected if any pairs of its points can be connected by a rational curve, see Kolla\'r $[17]$ for an elementary introduction to rationally connected varieties). Boyer, Hurtubise and Milgram $[1]$ proved that the stability theorems hold for principal almost solv-varieties, a subclass of rationally connected varieties on which a complex solvable linear group acts with a free dense open orbit.

In Symplectic setting, Cohen, Jones and Segal $[3]$ began to investigate what condition on a closed, connected, integral symplectic manifold $(X,\omega)$ ensures that the stability theorems hold. They consider the limit $\underset \longrightarrow  {\lim } \Hol_{x_0}(\C\P^1, X)$, which is viewed as a kind of stabilization of $\Hol_{x_0}(\C\P^1, X)$ under certain gluing operations. They conjecture that the stability theorem hold if and only if the evaluation map $E:\map {\underset \longrightarrow {\lim } \Hol_{x_0} (\C\P^1 , X)}. X.$
is a quasifibration.

All the above stability theorems have considered pairs of spaces at the same degree. A natural question is to ask that what is the connection among the spaces with varying degree. There are two reasons we want to ask such a question. One is that it is very natural in logic, and the other one is that it is easy to find a counterpart of the space $\Hol (\C\P^1, X)$ in algebraic geometry ( just $\Hom (\P^1, X)$ ) but $F(\P^1,X)$ makes no sense in general. So we have to consider the relations among the spaces with varying degree if we want to find some sort of stability property in algebraic geometry. Notice that homotopy in general also make no sense in algebraic geometry so that it should be replaced by something else. We have not found a good stability property for algebraic varieties over any algebraically closed field yet.

In this paper we will work over complex field $\C$ so that the concepts of continuous map and homotopy make sense.

Segal $[23]$ constructed a gluing map between $\Hol_{x_0}^{n} (\C\P^1, \C\P^1)$ and 
 $\Hol_{x_0}^{n+1}(\C\P^1,$

\noindent $\C\P^1)$ and showed that this gluing map induced isomorphisms 
 among the lower homology groups, i.e. he showed that $H_{*}(\Hol_{x_0}^{n} (\C\P^1, \C\P^1))$ is independent of $n$ when $n$ is large. I am not sure whether the similar stability theorem for homotopy groups was infered implicitly in his paper.

In general it is not so easy to find gluing maps between $\Hol_{x_0}^{*} (\C\P^1, X)$. However we find that the finite covering maps from $\C\P^1$ into itself provide a natural connection among $\Hol_{x_0}^{*} (\C\P^1, X)$, specifically we can define an injective map $g_{n,nd}: \map {\Hol_{x_0}^{n} (\C\P^1, X)}. {\Hol_{x_0}^{nd} (\C\P^1, X)}.$ which sends any $f \in \Hol_{x_0}^{n} (\C\P^1, X)$ to $g_{n,nd}(f)(z)= f(z^d), z\in \C\P^1$.

Now we can define the stability property for a variety.  

\begin{defn} We say that $X$ is a variety with the stability property if
and only if the injective map 
$$
g_{n,nd}:\map {\Hol_{x_0}^{n} (\P^1, X)}.{\Hol_{x_0}^{nd} (\P^1,X)}.
$$ 
is a rational homotopy equivalence up to dimension $k_n$ for any $x_0\in X$, where $k_n$ is a natural number dependent on $n$, $x_0$ and $X$, and $\underset {k\longrightarrow\infty} {\lim}
n_k=\infty$. 
\end{defn}

\begin{rem} Unfortunately the above injection does not induce homotopy 
equivalence among $\Hol_{x_0}^{*} (\P^1, X)$ even in the case $X=\C\P^1$, so we have to use rational homotopy instead of homotopy in the definition.
\end{rem} 

We wonder whether the following is true.

\begin{enumerate}

\item Is the stability property birationally invariant? That is if X is
birational to Y and X is a variety with stability property, so is Y?

\item Is it true that all rationally connected varieties are with stability property?
\end{enumerate}

Slight generally, if $X$ is a nilpotent complex algebraic variety with $\pi_2(X)$ a free abelian group of rank $r$ we can assign the homotopy class of any map $f \in \Hol (\C\P^1, X)$ a multiple degree $\bf{n}$. We can define the injective map $g_{{\bf n},  d {\bf n}}$ and 
stability in a similar way.

We conjecture that the following is true,
\begin{conj} Let $X$ be a variety as above. Then The map $g_{{\bf n},  d {\bf n}}$ induces a rational homotopy equivalence up to dimension $k_{\bf n}$, where $k_{\bf n}$ grows with $\bf{n}$, as $\bf{n}$ moves to infinity in a suitable positive cones.
\end{conj}

In this paper we will prove the above conjecture under some conditions:
\begin{thm} Let $X$ be a variety as above. Assume that the inclusion map 
 $\Hol_{x_0}^{{\bf n}} (\C\P^1, X) \subset F_{x_0}^{{\bf n}} (\C\P^1, X)$ induces rational homotopy equivalence up to dimension $k_{\bf n}$, where $k_{\bf n}$ grows with $\bf{n}$, as $\bf{n}$ moves to infinity in a suitable positive cones. Then $g_{{\bf n}, d {\bf n}}$ will induce rational homotopy equivalence up to dimension $\min \{k_{\bf n},k_{d{\bf n}}\}$.
\end{thm}
\begin{rem} The assumption that the inclusion map $\Hol_{x_0}^{{\bf n}} (\C\P^1, X) \subset F_{x_0}^{{\bf n}} (\C\P^1,$

\noindent $ X)$ induces rational homotopy equivalence up to dimension $k_n$ is not necessary. The theorem probably hold under much weaker conditions. Right now we do not know what can be used to replace this assumption. We also can not find a direct proof that $g_{{\bf n}, d {\bf n}}$ induces a rational homotopy equivalence up to dimension $\min \{k_{\bf n},k_{d{\bf n}}\}$. Here we will give an indirect proof of this by transfering the problem from the space of holomorphic maps into the space of continuous maps.

However, in a lot of cases (see the beginning of the introduction) we have known that the above assumption is true. So our theorem exactly gives a connection among the spaces with varing degrees.

Observe that Cohen, Jones and Segal $[3]$ defined some kind of stabilization of $\Hol_{x_0}(\C\P^1, X)$ under certain gluing operations. We can also define a stabilization of $\Hol_{x_0}(\C\P^1, X)$ by using the maps $g_{{\bf n}, d {\bf n}}$. The set $\{{\bf n}\}$ forms a partial order set. There is a morphism from ${\bf n}$ to ${\bf m}$ if and only if ${\bf m}=d {\bf n}$. This makes us to define the colimit $\underset \longrightarrow {\lim } \Hol_{x_0} (\P^1, X)$ over the directed system $\{{\bf n}\}$. We have not checked the connection between these two stabilizations of $\Hol_{x_0}(\C\P^1, X)$ yet. 
\end{rem}

The idea of the proof of Theorem $1.4$ is very simple. Look at the following commutative diagram,

$$
\begin{diagram}
 \Hol_{x_0}^{{\bf n}} (\C\P^1, X) &        & \rTo^i & & F_{x_0}^{{\bf n}} (\C\P^1, X) \\
\dTo^{g_{{\bf n}, d{\bf n}}} &       &   &     & \dTo_{g_{{\bf n}, d{\bf n}}} \\
\Hol_{x_0}^{d{\bf n}} (\C\P^1, X)          &  & \rTo^i &  & F_{x_0}^{d{\bf n}} (\C\P^1, X)       \\
\end{diagram}
$$

The horential arrows induce rational homotopy equivalence up to dimension $k_{\bf n}$ and $k_{d{\bf n}}$ respectively. If we can prove that the right arrow induces rational homotopy equivalence, then the left arrow induces rational homotopy equivalence up to dimension $\min \{k_{\bf n},k_{d{\bf n}}\}$.

We first show that the map $g_{{\bf n}, d{\bf n}}: \map F^{{\bf n}} (\C\P^1, X). F^{d{\bf n}} (\C\P^1, X).$ induces rational homotopy equivalence for free maps by using the Sullivan-Haefliger model in rational homotopy setting. The base map space $F_{x_0}^{*} (\C\P^1, X)$ is exactly the fiber of the evaluation map $E: \map F^{*} (\C\P^1, X). X.$ which sends $f$ to $f(\infty) \in X$ so that we can use the long exact sequence of the fibration to pass the rational homotopy equivalence from the space of free maps to the space of based maps. 

\noindent {\bf Acknowledgements} It is a great pleasure to express my appreciation to my advisor James M$^{c}$Kernan for proposing this problem. I am grateful to Yves F{\'e}lix and Samuel Bruce Smith for answering me some questions and pointing me out some important references on rational homotopy theory.

\section{Space of maps}

For any pair of topological spaces $X$ and $Y$, let $F(Y,X)$ be the space of all maps of $Y$ into $X$. In general, $F(Y,X)$ is a disconnected space, so for any map $f:\map Y.X.$, we let 
$F(Y,X,f)\subset F(Y,X)$ denote the path component that contains $f$. A fundamental problem is to classify the components of $F(Y,X)$ up to homotopy type. A lot of work $[4],[11],[12],[20],[21], [27]$ concerned the case that $X=S^n$, $\R\P^n$ and $\C\P^n$ has been done. A little weaker but still very interesting problem is to classify the components up to rational homotopy type. R. Thom $[28]$ studied the homotopy type of $F(Y,X,f)$ and computed explicitly the cohomology of $F(Y,X,f)$ when $X$ is a product of Eilenberg-Mac Lane spaces. Later on, following ideas of Sullivan, A. Haefliger $[10]$ gave the rational minimal model of the space of sections of a nilpotent bundle. This model has been extensively studied by K. Shibata $[24]$ and M. Vigu\'e-Poirrier $[29]$. Using a method of Haefliger's $[10]$, M$\varnothing$ller and Raussen $[22]$ proved that all components of $F(Y,S^n)$ are rationally homotopy equivalent if $n$ is an odd integer; and represent exact two rational homotopy types if $n$ is an even integer and $\tilde{H}^{0}(Y,\Q)=H^{\geq 2n-1}(Y,\Q)=0$. Their technique works well whenever the target space $Y$ is a $2$-stage Postnikov tower. As an illustration of this, they proved similar result for the space of maps into complex (and quaternionic) projective $n$-space.  In $[5]$ Y. F\'elix computed the rational category of the space of sections of a nilpotent bundle. S. B. Smith $[25], [26]$ also investigated the rational homotopy types of the components of some function spaces. For rational homotopy theory we refer to $[6]$.

Naturally we can identify $F(Y,X,f)$ with $\Gamma_{f'}$, the space
of sections of the trivial fibration $p_1: \map Y \times X. Y.$ homotopic to the section $f'
= (id, f): \map Y. Y \times X. $. To prove that the injection $\map F(Y, X, f).  F(Y, X, f \circ g).$ induces rational homotopy equivalence is equivalent to show that so does for $\map \Gamma_{f'}. \Gamma_{(f \circ g)'}.$, where $(f \circ g)'= (id, f \circ g)$ and $g(z)=z^d, z\in \C\P^1=S^2, d\in \N$.

Now let us recall the Sullivan-Haefliger model.

\section {The Sullivan-Haefliger model}

Let $X$ be a nilpotent space and $Y$ a simply connected space with a minimal model which is finite dimensional in each degree and $\dim H^*(Y)$ is finite. Let $p_1: \map {Y \times X}. Y.$ be the trival fibration. Denote $(A,d_A)$ be a model for $Y$ and $(\land V, d)$ be the minimal model for $X$. Let $f'$ be a section. It gives a morphism $\sigma_{f'}:\map {A \otimes \land V}. A.$ which is the identity on $A$. Let $\psi$ be the $A$-algebra automorphism of $A \otimes \land V$ mapping $1\otimes v$ on $1\otimes v -\sigma_{f'}(v)$; define on $A \otimes \land V$
a new differential $d'$ such that ${(d_A \otimes d)}\circ \psi= \psi \circ d'$. Then $\psi$ is a DG-automorphism and $\sigma_{f'} \circ \psi$ maps $\land V$ on zero. Denote the graded dual vector space of $A$ by $A^\vee$:
$$(A^\vee)_n= \Hom(A^n, \Q).$$
Let $(a_i)_{i \in I}$ be a basis of $A$. Then $A^\vee$ is naturally equipped with the dual basis $a_i^*$ such that $<a_i^*;a_j>=\delta_{ij}$.

We now look at the map of algebras defined by:
$$\varepsilon : \map {A \otimes \land V}. {A \otimes (A^\vee \otimes \land V)}.: \varepsilon (v) = \underset {i \in I} {\sum} {a_i \otimes (a_i^* \otimes v)}; \varepsilon (a)=a, a \in A$$

In $[10]$, Haefliger shows how to put a uniquely defined differential $d_A \otimes \delta$ on 
$A \otimes (A^\vee \otimes \land V)$ in such a way that $\varepsilon$ becomes a morphism of commutative differential graded algebras. Let $W$ be the quotient of $A^\vee \otimes \land V$ by the subspace of elements of degree $<0$, and by the subspace formed by the $\delta$-cocycles in degree $0$. Haefliger showed that $(\land W, \delta)$ is a model of the space $\Gamma_{f'}$. It has the following universal property. Let $D$ be a DG-algebra such that $D^0=\Q$, and let $f: \map {A \otimes \land V}. {A \otimes D}.$ be a morphism of augmented DG-algebras over $A$. Then there is a unique $\varphi: \map {\land W}. D.$ such that the diagram
$$
\begin{diagram}
  A \otimes \land V      &                                      & \rTo^f &       &  A \otimes D   \\                                                        
          & \rdTo_\varepsilon   &               & \ruTo_{1 \otimes \varphi}&          \\
          &                                      & A \otimes \land W             &                                        &          \\
\end{diagram}
$$
commutes.

Let $(A \otimes \land V, d_i), i=1,2$ be differential graded commutative algebras. Let $h:\map (A \otimes \land V, d_1). (A \otimes \land V, d_2).$ be an isomorphism which restricts to an automorphism on $A$. According to Haefliger's theorem, there are unique difined differential $\delta_i, i=1,2$ such that $\varepsilon_i, i=1, 2$ become morphisms of commutative differential graded algebras. Let $W_i, i=1, 2$ be constructed as above.

We have the following theorem:
\begin{thm} Notations as above. Then $(\land W_1, \delta_1)$ is isomorphic to 
$(\land W_2, \delta_2)$.
\end{thm}
\begin{proof}
We will apply the universal property of the Sullivan-Haefliger model to prove our theorem. Let us look at the following diagram
$$
\begin{diagram}
 (A \otimes \land V, d_1) & & \rTo^{\varepsilon_1} && (A \otimes \land W_1, d_A \otimes \delta_1) & &\rTo^{h|_{A} \otimes 1} && (A \otimes \land W_1, d_A \otimes \delta_1) \\
\dTo^{h} & &   & &   &   &    & &  \\
(A \otimes \land V, d_2)& &\rTo^{\varepsilon_2}& & (A \otimes \land W_2, d_A \otimes \delta_2) & &\rTo^{{h|_{A}}^{-1} \otimes 1}&  & (A \otimes \land W_2, d_A \otimes \delta_2) \\
\end{diagram}
$$

Obviously the map $(h|_{A} \otimes 1) \circ {\varepsilon_1} \circ (h^{-1}): (A \otimes \land V, d_2) \longrightarrow (A_1 \otimes \land V, d_1) \longrightarrow (A \otimes \land W_1, d_A \otimes \delta_1) \longrightarrow (A \otimes \land W_1, d_A \otimes \delta_1)$ is a morphism of augmented DG-algebras over $A$ and $\land W_1$ is a DG-algebra with $W_1^0=\Q$. Apply the universal property, there is a unique morphism $\varphi_{21}:\map {\land W_2}. {\land W_1}.$ such that $(1 \otimes {\varphi_{21}}) \circ {\varepsilon_2} = (h|_A \otimes 1) \circ {\varepsilon_1} \circ h^{-1}$. Similarly we can prove that there exist a unique morphism  $\varphi_{12}:\map {\land W_1}. {\land W_2}.$ such that $(1 \otimes \varphi_{12}) \circ {\varepsilon_1} =({h|_A}^{-1} \otimes 1) \circ {\varepsilon_2} \circ h: (A \otimes \land V, d_1) \longrightarrow (A_1 \otimes \land V, d_2) \longrightarrow (A \otimes \land W_2, d_A \otimes \delta_2) \longrightarrow (A \otimes \land W_2, d_A \otimes \delta_2)$. Now we apply the same trick as above to show that $\varphi_{21} \circ \varphi_{12} =1_{\land W_1}$ and $\varphi_{12} \circ \varphi_{21} =1_{\land W_2}$.

Since $(1 \otimes \varphi_{12}) \circ {\varepsilon_1} =({h|_A}^{-1} \otimes 1) \circ {\varepsilon_2} \circ h$, we have 
$(h|_A \otimes {\varphi_{21} \circ \varphi_{12}}) \circ {\varepsilon_1}= (1 \otimes \varphi_{21}) \circ (h|_A \otimes 1) \circ (1 \otimes \varphi_{12}) \circ {\varepsilon_1}= (1 \otimes \varphi_{21}) \circ (h|_A \otimes 1) \circ ({h|_A}^{-1} \otimes 1) \circ {\varepsilon_2} \circ h= (1 \otimes \varphi_{21}) \circ {\varepsilon_2} \circ h= (h|_A \otimes 1) \circ {\varepsilon_1} \circ h^{-1} \circ h = (h|_A \otimes 1) \circ {\varepsilon_1}$. So $(1 \otimes {\varphi_{21} \circ \varphi_{12}}) \circ {\varepsilon_1} = ({h|_A}^{-1} \otimes 1) \circ (h|_A \otimes {\varphi_{21} \circ \varphi_{12}}) \circ {\varepsilon_1}=({h|_A}^{-1} \otimes 1) \circ (h|_A \otimes 1) \circ {\varepsilon_1}= \varepsilon_1$

So there are two morphisms from $\land W_1$ to itself such that the following diagram
$$
\begin{diagram}
  A \otimes \land V &&  \rTo^{\varepsilon_1} && A \otimes \land W_1 \\
       &\rdTo_{\varepsilon_1}&       &  \ruTo_{1 \otimes (\varphi_{21} \circ \varphi_{12})}\ruTo^{1 \otimes 1} &          \\
          & &     A \otimes \land W_1   & &         \\
\end{diagram}
$$
commutes. By the uniqueness of such morphism, we have $\varphi_{21} \circ \varphi_{12}=1_{\land W_1}$. Similarly we can prove $\varphi_{12} \circ \varphi_{21}=1_{\land W_2}$. This completes our proof.

\end{proof}

Applying the above theorem we can prove
\begin{thm} The map $g_{{\bf n}, d{\bf n}}: \map F^{{\bf n}} (\C\P^1, X). F^{d{\bf n}} (\C\P^1, X).$ is a rational homotopy equivalence.
\end{thm}
\begin{proof}
Applying Theorem $3.1$, We have $Y=\C\P^1=S^2$. Let $(A,d_A)$ be a minimal model for $S^2$. Let $g(z)=z^d, z\in \C\P^1=S^2, d\in \N$. It induces an automorphism of $(A,d_A)$, denote as $\sigma_g$. For the section $f'$ and $(f \circ g)'$, the induced morphisms on $(A \otimes \land V, d_A \otimes d)$ are $\sigma_{f'}$ and $\sigma_{(f \circ g)'}$ respectively. Just as in Haefliger $[10]$ we have two differential graded algebra $(A \otimes \land V, d_i), i=1, 2$, where $d_i, i=1, 2$ are definded by ${(d_A \otimes d)}\circ \psi_i= \psi_i \circ d_i, i=1, 2$ and $\psi_1$ and $\psi_2$ mapping $a \in A$ to itself and $1\otimes v$ on $1\otimes v -\sigma_{f'}(v)$ and on $1\otimes v -\sigma_{(f \circ g)'}(v)$, respectively.

Now look at the following diagram,
$$
\begin{diagram}
(A \otimes \land V, d_A \otimes d) & &\rTo^{\psi_1}& & (A \otimes \land V, d_1)\\
\dTo^{\sigma_g \otimes 1} & &                            &  & \dTo_{\sigma_g \otimes 1}\\
(A \otimes \land V, d_A \otimes d) && \rTo^{\psi_2} && (A \otimes \land V, d_2)\\
\end{diagram}
$$

We claim that the above is a commutative diagram. To show this, let $a \otimes v$ be any
element in $A \otimes \land V$. $\psi_1 (a \otimes v)=a \otimes v -a \sigma_{f'}(v)$. So $(\sigma_g \otimes 1) \circ \psi_1 (a \otimes v)= \sigma_{g}(a) \otimes v-\sigma_{g}(a)\sigma_{g}(\sigma_{f'}(v))$. On the other hand $\psi_{2} \circ (\sigma_g \otimes 1) (a \otimes v)=\psi_{2}(\sigma_{g}(a) \otimes v)= \sigma_{g}(a) \otimes v- \sigma_{g}(a) \sigma_{(f\circ g)'}(v)$. If we can prove that $\sigma_{g}(\sigma_{f'}(v))= \sigma_{(f\circ g)'}(v)$, then our claim follows.
Let $p_2: \map {S^2 \times X}. X.$ be the second projection. It induce a map $\sigma_{p_2}: \map {\land V}. {A \otimes \land V}.$ which maps $v$ on $1 \otimes v$. Since $f=p_2 \circ f'$, so $\sigma_f=\sigma_{f'} \circ \sigma_{p_2}: \map \land V. A.$, so $\sigma_{f}(v)=\sigma_{f'}(v)$. Similarly we can prove that $\sigma_{f \circ g}(v)=\sigma_{(f\circ g)'}(v)$. But $\sigma_{f \circ g}(v)=\sigma_{g}(\sigma_{f}(v))$. So $\sigma_{g}(\sigma_{f'}(v))= \sigma_{(f\circ g)'}(v)$. Since the left vertical arrow and the two horential arrows are isomorphisms of differential graded commutative algebras, so the right vertical arrows gives an isomorphism between $(A \otimes \land V, d_1)$ and $(A \otimes \land V, d_2)$. Obviously the map $(\sigma_g \otimes 1)|_A$ is an automorphism of $A$. So apply theorem 3.1 we can prove the Sullivan-Haefliger models for the space $\Gamma_{f'}$ and $\Gamma_{(f \circ g)'}$ are isomorphic. Therefore these two spaces have the same rational homotopy type.
\end{proof}

Now Theorem $1.4$ follows easily by the idea giving at the end of the introduction. We only need to consider the long exact homotopy sequences of the fibrations $F^{\bf n}_{x_0} (\C\P^1, X) \longrightarrow F^{\bf n}(\C\P^1, X) \longrightarrow X$ and  $F^{d{\bf n}}_{x_0} (\C\P^1, X) \longrightarrow F^{d{\bf n}}(\C\P^1, X) $

\noindent $\longrightarrow X$ and apply the Five-Lemma.

It is also interesting to note that there are only at most two rational homotopy type of the components of $F(\C\P^1, X)$ if $X$ is simply connected and $\pi_2(X)=\Z$.

\bigskip
\begin{center}  
{\bf References}
\end{center}
\bigskip

\begin{enumerate}
\item C. P. Boyer, J. C. Hurtubise, and R. J. Milgram, Stability theorems for spaces of rational curves, Internat. J. Math., {\bf 12}, 2001, No. 2, 223--262.

\item C. P. Boyer, B. M. Mann, J. C. Hurtubise and R. J. Milgram, The topology of the space of rational maps into generalized flag manifolds, Acta Math., {\bf 173}, 1994, No. 1, 61--101.

\item R. L. Cohen, J. D. S. Jones and G. B. Segal, Stability for holomorphic spheres and {M}orse theory, in Geometry and topology: Aarhus (1998), Contemp. Math., {\bf 258}, 87--106, Amer. Math. Soc., Providence, RI, 2000.

\item M. C. Crabb and W. A. Sutherland, Function spaces and {H}urwitz-{R}adon numbers, Math. Scand., {\bf 55}, 1984, No. 1, 67--90.

\item Y. F{\'e}lix, Rational category of the space of sections of a nilpotent bundle, Comment. Math. Helv., {\bf 65}, 1990, No. 4, 615--622.

\item Y. F{\'e}lix, S. Halperin and J.C. Thomas, Rational homotopy theory, Graduate Texts in Mathematics, {\bf 205}, Springer-Verlag, New York, 
 2001.

\item J. Gravesen, On the topology of spaces of holomorphic maps, Acta Math., {\bf 162}, 1989, No. 3-4, 247--286.

\item M. A. Guest, Topology of the space of absolute minima of the energy functional, Amer. J. Math., {\bf 106}, 1984, No. 1, 21--42.

\item M. A. Guest, The topology of the space of rational curves on a toric variety, Acta Math., {\bf 174}, 1995, No. 1, 119--145.
 
\item A. Haefliger, Rational homotopy of the space of sections of a nilpotent bundle, Trans. Amer. Math. Soc., {\bf 273}, 1982, No. 2, 609--620.

\item V. L. Hansen, The homotopy problem for the components in the space of maps on the {$n$}-sphere, Quart. J. Math. Oxford Ser. (2), {\bf 25}, 1974, 313--321.

\item V. L. Hansen, On spaces of maps of {$n$}-manifolds into the {$n$}-sphere, Trans. Amer. Math. Soc., {\bf 265}, 1981, No. 1, 273--281.

\item J. C. Hurtubise, Holomorphic maps of a {R}iemann surface into a flag manifold, J. Differential Geom. {\bf 43}, 1996, No. 1, 99--118.

\item J. C. Hurtubise, Moduli spaces and particle spaces, in Gauge theory and symplectic geometry (Montreal, PQ, 1995), Proc. NATO Adv. Sci. Inst. Ser. C Math. Phys. Sci., {\bf 488}, 113--146, Kluwer Acad. Publ., Dordrecht ,1997.

\item F. Kirwan, On spaces of maps from {R}iemann surfaces to {G}rassmannians and applications to the cohomology of moduli of vector bundles, Ark. Mat., {\bf 24}, 1986, No. 2, 221--275.

\item J. Koll\'ar, Rational curves on algebraic varieties, Ergebnisse der Mathematik und ihrer Grenzgebiete. 3. Folge. A Series of Modern Surveys in Mathematics [Results Mathematics and Related Areas. 3rd Series. A Series of Modern Surveys in Mathematics], 32, Springer-Verlag, 1996, viii+320.

\item J. Koll\'ar, Which are the simplest algebraic varieties?, Bull. Amer. Math. Soc. (N.S.), {\bf 38}, 2001, No. 4, 409--433.

\item B. M. Mann and R. J. Milgram, Some spaces of holomorphic maps to complex {G}rassmann
manifolds, J. Differential Geom., {\bf 33}, 1991, No. 2, 301--324.

\item B. M. Mann and R. J. Milgram, On the moduli space of {${\rm SU}(n)$} mono-
poles and
holomorphic maps to flag manifolds, J. Differential Geom., {\bf 38}, 1993, No. 1, 39--103.

\item J. M. M{\o}ller, On spaces of maps between complex projective spaces, Proc. Amer. Math. Soc., {\bf 91}, 1984, No. 3, 471--476.

\item J. M. M{\o}ller, On the homology of spaces of sections of complex projective bundles,
Pacific J. Math., {\bf 116}, 1985, No. 1, 143--154.

\item J. M. M{\o}ller and M. Raussen, Rational homotopy of spaces of maps into spheres and complex projective spaces, Trans. Amer. Math. Soc., {\bf 292}, 1985, No. 2, 721--732.

\item G. Segal, The topology of spaces of rational functions, Acta Math., {\bf 143}, (1979), No. 1-2, 39--72.

\item K. Shibata, On {H}aefliger's model for the {G}el\' fand-{F}uks cohomology, Japan. J. Math. (N.S.), {\bf 7}, 1981, No. 2, 379--415.

\item S. B. Smith, Rational homotopy of the space of self-maps of complexes with finitely many homotopy groups, Trans. Amer. Math. Soc., {\bf 342}, 1994, No. 2, 895--915.

\item S. B. Smith, Rational classification of simple function space components for flag manifolds, Canad. J. Math., {\bf 49}, 1997, No. 4, 855--864.

\item W. A. Sutherland, Path-components of function spaces, Quart. J. Math. Oxford Ser. (2), {\bf 34}, 1983, No. 134, 223--233.

\item R. Thom, L'homologie des espaces fonctionnels, in Colloque de topologie 
alg\'ebrique, Louvain, 1956, 29--39, Georges Thone, Li\`ege, 1957.

\item M. Vigu{\'e}-Poirrier, Sur l'homotopie rationnelle des espaces fonctionnels, 

\noindent Manuscripta Math., {\bf 56}, 1986, No. 2, 177--191.
\end{enumerate}
\medskip
\noindent Department of Mathematics, Syracuse University, NY 13244.
jilin@syr.edu

\end{document}